\documentclass{amsart}

\usepackage{amsmath}
\usepackage{amsfonts}
\usepackage[utf8]{inputenc}
\usepackage[numbers]{natbib}
\usepackage{graphicx}
\usepackage{mathtools}
\usepackage{amssymb}
\usepackage{amsthm}
\usepackage{tikz-cd} 
\usepackage{enumerate}
\usepackage[english]{babel}
\usepackage{hyperref}
\usepackage{tikz}
\usetikzlibrary{positioning} 
\usepackage{nomencl}

\theoremstyle{definition}
\newtheorem{cor}{Corollary}[section] 
\newtheorem{prop}[cor]{Proposition}
\newtheorem{theorem}[cor]{Theorem}

\newtheorem{conj}[cor]{Conjecture}
\newtheorem{defi}[cor]{Definition}

\newtheorem{lemma}[cor]{Lemma}

\numberwithin{equation}{section}

\DeclareMathOperator{\aut}{Aut}

\DeclareMathOperator{\su}{SU}

\DeclareMathOperator{\slg}{SL}

\DeclareMathOperator{\gu}{GU}

\DeclareMathOperator{\iso}{Iso}

\DeclareMathOperator{\syl}{Syl}
\DeclareMathOperator{\gl}{GL}
\DeclareMathOperator{\ob}{Ob}
\DeclareMathOperator{\spin}{Spin}

\DeclareMathOperator{\inj}{Inj}

\DeclareMathOperator{\out}{Out}
\DeclareMathOperator{\sym}{Sym}

\DeclareMathOperator{\id}{id}

\DeclareMathOperator{\homo}{Hom}
\DeclareMathOperator{\inn}{Inn}

\DeclareMathOperator{\br}{Br}

\newcommand{\bigslant}[2]{{\raisebox{.2em}{$#1$}\left/\raisebox{-.2em}{$#2$}\right.}}

\begin{document}


\title{Block-exoticity of a family of exotic fusion systems}


\author{Patrick Serwene}
\address{City University London, Northampton Square, London, EC1V 0HB, UK}
\email{patrick.serwene@city.ac.uk}





\begin{abstract}
We prove that each exotic fusion system $\mathcal F$ on a Sylow $p$-subgroup of $G_2(p)$ for an odd prime $p$ with $\mathcal O_p(\mathcal F)=1$ is block-exotic. This gives evidence for the conjecture that each exotic fusion system is block-exotic. We prove two reduction theorems for block-realisable fusion systems.
\end{abstract}






\maketitle

\section{Introduction}

Let $p$ be a prime, $P$ a finite $p$-group and $\mathcal F$ a fusion system on $P$. By fusion system we always mean saturated fusion system. Recall that $\mathcal F$ is said to be realisable if $\mathcal F=\mathcal F_P(G)$ for a finite group $G$ and $P \in \syl_p(G)$, otherwise $\mathcal F$ is said to be exotic. Also recall that if $\mathcal F=\mathcal F_{(P,e_P)}(G,b)$ for a finite group $G$ having a $p$-block $b$ with maximal $b$-Brauer pair $(P,e_P)$, $\mathcal F$ is said to be block-realisable, otherwise $\mathcal F$ is said to be block-exotic. See Section $2$ for details.\\
The following fact is a consequence of Brauer's Third Main Theorem (see \cite[Theorem 7.1]{radha}): If $G$ is a finite group and $b$ is the principal $p$-block of $kG$, i.e. the block corresponding to the trivial character, with maximal $b$-Brauer pair $(P,e_P)$, then $P \in \syl_p(G)$ and $\mathcal{F}_{(P,e_P)}(G,b)=\mathcal{F}_P(G)$. In particular, any realisable fusion system is block-realisable. The converse is still an open problem. However, we have the following

\begin{conj}
\label{seitan}
\textit{If $\mathcal F$ is an exotic fusion system, then $\mathcal F$ is block-exotic.}
\end{conj}

There have been only two families of exotic fusion systems for which block-exoticity has been proven. The first one being the Solomon systems defined on a Sylow $2$-subgroup of $\spin(q)$. These are conjectured to be the only exotic fusion systems on $2$-groups. The block-exoticity for $q=3$ was proven in \cite{radha} and generalised to all odd prime powers $q$ in \cite[Theorem 9.34]{david}. The second example consists of the Ruiz--Viruel systems, which are defined on the extra special $7$-group of order $7^3$ and for which block-exoticity has been proven in \cite{ks}.\\
In this paper, we provide further evidence for Conjecture \ref{seitan}. Our first main result is the following: 

\begin{theorem}
\label{oh hai mark}
\textit{Conjecture $\ref{seitan}$ is true for all fusion systems $\mathcal F$ on a Sylow $p$-subgroup of $G_2(p)$ for odd $p$ with $\mathcal O_p(\mathcal F)=1$}.
\end{theorem}

The fusion systems on a Sylow $p$-subgroup of $G_2(p)$ for odd $p$ and $\mathcal O_p(\mathcal F)=1$ have been classified by Parker and Semeraro in \cite{parkersem} and thus we refer to them as Parker--Semeraro systems. For $p \neq 7$, all Parker--Semeraro systems are realised by finite groups, whereas for $p=7$, there are $29$ Parker--Semeraro systems of which $27$  are exotic. In this paper, we prove block-exoticity of these.\\
The relevance of the Parker--Semeraro systems stems from the fact that one wants to classify all fusion systems over maximal unipotent subgroups of finite groups of Lie type of rank 2. Note that the Solomon systems belong to this class of fusion systems as well. Furthermore, another important factor is that $7$ is a good prime and thus many results for groups of Lie type will be applicable.\\
The proof of Theorem \ref{oh hai mark} uses the classification of finite simple groups. A key step is a reduction theorem which we state after the following definition:

\begin{defi}
\textit{Let $\mathcal F$ be a fusion system on a finite $p$-group $P$. If $P$ has no non-trivial proper strongly $\mathcal F$-closed subgroups, we call $\mathcal F$ reduction simple}. 
\end{defi}

Note that if $\mathcal F=\mathcal F_{(P,e_P)}(G,b)$, $b$ is said to be an $\mathcal F$-block. The reduction we apply to prove Theorem \ref{oh hai mark} takes the following form:

\begin{theorem}
\label{cheeze}
\textit{Let $\mathcal F$ be a reduction simple fusion system on a non-abelian $p$-group $P$. If $\mathcal F$ is block-realisable, then there exists a finite group $G$ possessing an $\mathcal F$-block $b$ such that the following holds\\
$(a)$ $|G : Z(G)|$ is minimal among all groups $G$ having an $\mathcal F$-block,\\
$(b)$ if $H \unlhd G$ with $P \not\subseteq H$, then $H$ is a central $p'$-group. In particular $F(G)=Z(G)$,\\
$(c)$ the number of components of $G$ is bounded by the rank of $Z(P)$}.
\end{theorem}

Our final main result is another reduction theorem, which generalises \cite[Theorem 4.2]{ks}. We refer to Section 2 for terminology.

\begin{theorem}
\label{youwereneverreallyhere}
\textit{Let $\mathcal{F}$ be a fusion system on a non-abelian $p$-group $P$. Denote by $\mathcal F_H$ the subsystem of $\mathcal F$ corresponding to $H \leq \Gamma_{p'}(\mathcal F)$. Assume\\
$(a)$ $\mathcal O^{p'}(\mathcal F)$ is reduction simple,\\
$(b)$ if $\mathcal G$ is a fusion system on $P$ containing $\mathcal{O}^{p'}(\mathcal F)$, then $\mathcal G \subseteq \mathcal F$,\\
$(c)$ if $\mathcal G$ is a fusion system on $P$ such that $\mathcal G \unlhd \mathcal F_H$, then $\mathcal{O}^{p'}(\mathcal{F}) \subseteq \mathcal{G}$.\\
If there exists a finite group having an $\mathcal F_H$-block, then there exists a finite quasisimple group $L$ with $p'$-centre having an $\mathcal F_{H'}$-block for some $H' \leq \Gamma_{p'}(\mathcal F)$}.
\end{theorem}

In the next section, we recall the definition of fusion systems, several of their key properties and results. Section 3 will be concerned with reduction results and we prove Theorems \ref{cheeze} and \ref{youwereneverreallyhere}. Finally, Section 4 will be about the Parker--Semeraro systems and we prove Theorem \ref{oh hai mark}.

\section{Background and quoted results on fusion systems}
We begin by recalling the definition of a fusion system.

\begin{defi}
\label{jabloka}
\textit{Let $p$ be a prime and $P$ be a finite $p$-group. A saturated fusion system on $P$ is a category $\mathcal F$, s.t. $\ob(\mathcal F)$ is the set of all subgroups of $P$ and furthermore for all $Q, R \leq P$ we have:\\
$(i) \ \homo_P(Q,R) \subseteq \homo_{\mathcal F}(Q,R) \subseteq \inj(Q,R),$\\
$(ii)$ each homomorphism in $\mathcal F$ is the composition of an $\mathcal F$-isomorphism and an inclusion$,$\\
$(iii)$ each subgroup of $P$ is $\mathcal F$-conjugate to a subgroup which is fully automised and receptive in $\mathcal F$.}
\end{defi}

We need to recall some of the notation from this definition:

\begin{defi}
\textit{$(a)$ Let $\mathcal F$ be a fusion system on a $p$-group $P$. Two subgroups $Q,R \leq P$ are called $\mathcal F$-conjugate if they are isomorphic as objects of the category $\mathcal F.$\\
$(b)$ A subgroup $Q \leq P$ is called fully automised in $\mathcal F$ if $\aut_P(Q) \in \syl_p(\aut_{\mathcal F}(Q)).$\\
$(c)$ A subgroup $Q \leq P$ is called receptive in $\mathcal F$ if for each $R \leq P$ and each $\varphi \in \iso_{\mathcal{F}}(R,Q)$, $\varphi$ has an extension to the group $N_{\varphi}:=\{ g \in N_P(R) \mid$ $^{\varphi} c_{g} \in \aut_P(Q) \}.$}
\end{defi}

For convenience, we drop the term saturated, and mean saturated fusion system whenever we say fusion system in this paper. In the literature, fusion system means categories satisfying only axioms $(i)$ and $(ii)$ from Definition \ref{jabloka}. Furthermore, note that a category whose objects consist of the subgroups of some $p$-group $P$ is called category on $P$.

\begin{theorem}
\textit{Let $G$ be a finite group with $P \in \syl_p(G)$. We denote the category on $P$ with morphisms consisting of conjugation by elements in $G$ by $\mathcal{F}_P(G)$. Then $\mathcal F_P(G)$ is a fusion system on $P$}.
\end{theorem}

If a fusion system is of the form defined in the previous theorem, we call it realisable, otherwise we call it exotic.

\begin{defi}
\textit{
Let $\mathcal F$ be a fusion system on a $p$-group $P$ and $\mathcal E \subseteq \mathcal F$ be a subcategory of $\mathcal F$ which is a fusion system itself on some subgroup $P' \leq P$.\\
$(a)$ A subgroup $Q \leq P$ is called strongly $\mathcal F$-closed, if $\varphi(R) \subseteq Q$ for each $\varphi \in \homo_{\mathcal F}(R,P)$ and each $R \leq Q$.\\
$(b)$ If $P'$ is normal in $P$ and strongly $\mathcal F$-closed, $^\alpha \mathcal E=\mathcal E$ for each $\alpha \in \aut_{\mathcal F}(P')$ and for each $Q \leq P'$ and $\varphi \in \homo_{\mathcal F}(Q,P')$, there are $\alpha \in \aut_{\mathcal F}(P')$ and $\varphi_0 \in \homo_{\mathcal E}(Q, P')$ with $\varphi = \alpha \circ \varphi_0$, then $\mathcal E$ is called weakly normal in $\mathcal F$, denoted $\mathcal E \dot{\unlhd} \mathcal F$.\\
$(c)$ If $\mathcal E$ is weakly normal and in addition, we have that each $\alpha \in \aut_{\mathcal E}(P')$ has an extension $\overline{\alpha} \in \aut_{\mathcal F}(P'C_P(P'))$ with $[\overline{\alpha}, C_P(P')] \leq Z(P')$, then we call $\mathcal E$ normal in $\mathcal F$ and write $\mathcal E \unlhd \mathcal F$.\\
$(d)$ A fusion system is called simple if it does not contain any non-trivial proper normal subsystem}.
\end{defi}

The collection of weakly normal subsystems of $\mathcal F$ on $P$, ordered by inclusion, has a unique minimal element, see \cite[Theorem 7.7]{ako}. We denote this subsystem by $\mathcal{O}^{p'}(\mathcal F)$. Additionally, if $\mathcal F$ is a fusion system on $P$ and $M \subseteq \aut_{\mathcal F}(Q)$ for some $Q \leq P$, we denote by $\langle M \rangle$ the smallest (not necessarily saturated) fusion system on $P$ such that its morphisms contain $M$.

\begin{defi}
\textit{$(a)$ Let $\mathcal F$ be a fusion system on a finite $p$-group $P$ and $Q \leq P$. If $C_{P}(Q')=Z(Q')$ for each $Q' \leq P$ which is $\mathcal F$-conjugate to $Q$, then $Q$ is called $\mathcal F$-centric.\\
$(b)$ We say that a subsystem of $\mathcal F$ has index prime to $p$ $($or $p'$-index$)$ if it contains $\mathcal O^{p'}(\mathcal F)$.\\
$(c)$ We define the groups $\mathcal O^{p'}_{\ast}(\mathcal F):=\langle O^{p'}(\aut_{\mathcal F}(Q)) \mid Q \leq P \rangle $, $\aut^0_{\mathcal F}(\mathcal F):=\langle \alpha \in \aut_{\mathcal F}(P) \mid$ $\alpha \mid_{Q}$ $\in \homo_{\mathcal O^{p'}_{\ast}(\mathcal F)}(Q,P)$ for some $\mathcal F$-centric $Q \leq P\rangle$, and $\Gamma_{p'}(\mathcal F):=\aut_{\mathcal F}(P)/\aut^0_{\mathcal F}(P)$}.
\end{defi}

Note that the above definition makes sense since clearly $\mathcal O^{p'}_{\ast}(\mathcal F) \leq \mathcal F$ and also $\aut^0_{\mathcal F}(P) \unlhd \aut_{\mathcal F}(P)$.

It turns out that the group $\Gamma_{p'}(\mathcal F)$ carries a lot of information about $\mathcal F$:

\begin{theorem}
\label{mabel}
\cite[Theorem 7.7]{ako} \textit{Let $\mathcal F$ be a fusion system on a $p$-group $P$. There is a one-to-one-correspondence between the subsystems of $\mathcal F$ of index prime to $p$ and subgroups of $\Gamma_{p'}(\mathcal F)$. For some $H \leq \Gamma_{p'}(\mathcal F)$, we will refer to the corresponding subsystem by $\mathcal F_H$. Furthermore, this correspondence respects (weak) normality}.
\end{theorem}

For the fusion systems we encounter in this article, the group $\Gamma_{p'}(\mathcal F)$ will always be cyclic.

\begin{defi}
\textit{Let $\mathcal F$ be a fusion system on a $p$-group $P$. We call a subgroup $Q \leq P$ normal in $\mathcal F$, denoted $Q \unlhd \mathcal F$, if $Q \unlhd P$ and any morphism $\varphi \in \homo_{\mathcal F}(R,S)$ for $R,S \leq P$ has an extension $\overline \varphi \in \homo_{\mathcal F}(RQ,SQ)$ with $\overline \varphi(Q)=Q$. The largest subgroup of $P$ normal in $\mathcal F$ is denoted by $\mathcal{O}_p(\mathcal F)$}.
\end{defi}

\begin{lemma}
\label{pint}
\textit{Let $\mathcal F, \mathcal G$ be fusion systems on a finite $p$-group $P$ such that $\mathcal F \leq \mathcal G$ and let $Q \unlhd P$. If $Q$ is normal in $\mathcal G$, then $Q$ is normal in $\mathcal F$}.
\end{lemma}

\textit{Proof}. By \cite[Proposition 4.5]{ako}, this is equivalent to show that if $Q$ is contained in each $\mathcal F$-essential subgroup $R$ of $P$ and for each of these $R$, we have that $Q$ is $\aut_{\mathcal F}(R)$-invariant. Since $Q$ is normal in $\mathcal G$, it is strongly $\mathcal F$-closed. In particular,  $Q$ is $\aut_{\mathcal F}(R)$-invariant for all $R \leq P$ such that $Q \leq R$.\\
Now let $\varphi: R \rightarrow T$ be an $\mathcal G$-isomorphism. For the other property, we first claim $N_P(R) \cap Q \leq N^\mathcal{G}_\varphi$. Indeed, since $Q \unlhd \mathcal G$, $\varphi$ extends to $\overline{\varphi}: QR \rightarrow QT$. Now for $n \in N_P(R) \cap Q$, if we consider the map $c_n: R \rightarrow R$ then $^\varphi c_n = \varphi \circ c_n \circ \varphi^{-1}=c_{\overline{\varphi}(n)} \in \aut_P(T)$, which means $n \in N^\mathcal{G}_\varphi$.\\
Now let $R \leq P$ be $\mathcal F$-essential and $\beta \in \aut_{\mathcal F}(R)$ such that $N^\mathcal{F}_\beta = R$ (such an $\beta$ exists since $R$ is $\mathcal F$-essential). One easily verifies $N^\mathcal{F}_\beta=N^\mathcal{G}_\beta$. So, by the above, we have $N_P(R) \cap Q \leq N^\mathcal{F}_\beta = R$. Since $Q \unlhd P$, we have $RQ \leq P$. By general properties of $p$-groups either $RQ=R$ or $R < N_{RQ}(R)$. Since we have $N_{RQ}(R)=RN_Q(R)=R(N_P(R) \cap Q)=R$, we deduce $RQ=R$, so $Q \leq R$, which implies the claim. \hfill$\square$

\begin{lemma}
\label{ofscience}
\textit{If $\mathcal F, \mathcal G$ are fusion systems on a finite $p$-group $P$ with $\mathcal G \dot{\unlhd} \mathcal F$, then $O_p(\mathcal G)$ is normal in $\mathcal F$}.
\end{lemma}

\textit{Proof}. We have to check that each morphism $\varphi \in \homo_\mathcal{F}(Q,P)$ has an extension $\overline \varphi: QO_p(\mathcal G) \rightarrow P$ with $\overline \varphi(O_p(\mathcal G))=O_p(\mathcal G)$. Since $\mathcal G \dot{\unlhd} \mathcal F$, $\varphi$ can be written as $\varphi = \alpha \circ \beta$, where $\alpha \in \aut_\mathcal{F}(P)$ and $\beta \in \homo_\mathcal{G}(Q,P)$. By Alperin's Fusion Theorem, see \cite[Theorem 5.2]{markus}, we may assume $\beta$ is a bijection. So let $\tau $ be such that $\beta \circ \tau=\id_Q$. By definition of $O_p(\mathcal G)$, $\beta$ has an extension as desired. Thus, after replacing $\varphi$ with $\varphi \circ \tau$, we can assume $\varphi \in \aut_\mathcal{F}(P)$. It remains to show that $\varphi(O_p(\mathcal G))$ is normal in $\mathcal G$ since that implies $\varphi(O_p(\mathcal G)) \subseteq O_p(\mathcal G)$.\\
Let $\psi: R \rightarrow P$ be a morphism in $\mathcal{G}$. We need to show that $\psi$ extends to a morphism $\overline \psi: \varphi(O_p(\mathcal G)) R \rightarrow P$ with $\overline \psi (\varphi(O_p(\mathcal G)))=\varphi(O_p(\mathcal G))$. Now $\varphi^{-1} \psi \varphi: \varphi^{-1}(R) \rightarrow P$ is a morphism in $\mathcal G$ by weak normality. In particular, it has an extension $\gamma: \varphi^{-1}(R) O_p(\mathcal G) \rightarrow P$ with $\gamma(O_p(\mathcal G))=O_p(\mathcal G)$. Define $\pi:=\varphi \gamma \varphi^{-1}: R \varphi(O_p(\mathcal G)) \rightarrow P$. Then for $r \in R$, we have $\pi(r)=\varphi \gamma \varphi^{-1}(r)=\varphi \varphi^{-1} \psi \varphi \varphi^{-1}(r)=\psi(r)$ and for $n \in O_p(\mathcal G)$, we have $\pi(\varphi(n))=\varphi \gamma \varphi^{-1}(\varphi(n))=\varphi \gamma(n) \in \varphi(O_p(\mathcal G))$. This means we can take $\pi$ as $\overline \psi$. \hfill$\square$\\

Since fusion systems which do not allow many strongly closed subgroups play an important role, we define the following:

\begin{defi}
\textit{Let $\mathcal F$ be a fusion system on a $p$-group $P$. If $P$ does not have any non-trivial proper strongly $\mathcal F$-closed subgroups, we call $\mathcal F$ reduction simple}.
\end{defi}

The rest of this section will recollect some block-theoretic results. In particular, the next theorem shows how blocks of finite groups provide fusion systems.

\begin{defi}
\textit{Let $G$ be a finite group and $b$ a block of $kG$. A Brauer pair is a pair $(Q,f)$ where $Q \leq_p G$ and $f$ is a block of $kC_G(Q)$. We denote the set of blocks of $kC_G(Q)$ for some $p$-subgroup $Q$ of $G$ by $\mathcal B(Q)$}.
\end{defi}

Note that $G$ acts on the set of Brauer pairs by conjugation. Furthermore, Brauer pairs form a poset:

\begin{defi}
\textit{Let $(Q,f)$ and $(R,e)$ be Brauer pairs. Then\\
$(a) \ (Q,f) \unlhd (R,e)$ if $Q \unlhd R$, $f$ is $R$-stable and $\br_R(f)e=e,$\\
$(b) \ (Q,f) \leq (R,e)$ if $ Q \leq R$ and there exist Brauer pairs $(S_i,d_i), 1 \leq i \leq n$ such that $(Q,f) \unlhd (S_1, d_1) \unlhd (S_2, d_2) \unlhd \dots \unlhd (S_n,d_n) \unlhd (R,e)$}.
\end{defi}

It is actually the case that for a given Brauer pair $(R,e)$ and some $Q \leq R$, there exists a unique $f \in \mathcal{B}(Q)$ with $(Q,f) \leq (R,e)$, see \cite[Theorem 2.9]{radha}.

\begin{defi}
\textit{Let $b$ be a block of a finite group $G$.\\
$(a)$ A $b$-Brauer pair is a Brauer pair $(R,e)$ such that $(1,b) \leq (R,e)$, or equivalently it is a Brauer pair $(R,e)$ such that $\br_R(b)e=e.$\\
$(b)$ We denote the blocks $e$ of $kC_G(R)$ such that $(1,b) \leq (R,e)$ by $\mathcal B(R,b).$\\
$(c)$ A defect group of $b$ is a $p$-subgroup $P$ of $G$ maximal such that $\br_P(b) \neq 0$.}
\end{defi}

Note that the group $G$ acts by conjugation on the set of $b$-Brauer pairs. Furthermore, some $P \leq_p G$ is a defect group of $b$ if and only if there is a maximal pair $(P,e)$ such that $(1,b) \leq (P,e)$. We refer to such a pair as a maximal $b$-Brauer pair.

\begin{theorem}
\textit{Let $b$ be a block of $kG$ and $(P,e_P)$ be a maximal $b$-Brauer pair. For a subgroup $Q \leq P$, denote by $e_Q$ the unique block such that $(Q,e_Q) \leq (P,e_P)$. Denote the category on $P$ whose morphisms consist of all injective group homomorphisms $\varphi: Q \rightarrow R$ for which there is some $g \in G$ such that $\varphi(x)$ = $^gx$ for all $x \in Q$ and $^g(Q,e_Q) \leq (R,e_R)$ by $\mathcal F_{(P,e_P)}(G,b)$. Then $\mathcal F_{(P,e_P)}(G,b)$ is a fusion system on $P$}.
\end{theorem}

If a fusion system is of the form defined in the previous theorem, we call it block-realisable, otherwise we call it block-exotic.

Note that for several of the reduction theorems we need to introduce more general structures than block fusion systems, since some group theoretic properties are not captured by these: Assume $b$ is a block of $kG$ with maximal $b$-Brauer pair $(P,e)$ and $N \unlhd G$. If $c$ is a block of $kN$ covered by $b$, i.e. $bc \neq 0$, $P \cap N$ will be a defect group for $c$. However, in general $\mathcal F_{(P \cap N, e_{P \cap N})}(N,c)$ does not even need to be a subsystem of $\mathcal F_{(P,e_P)}(G,b)$. We introduce a generalised category to circumvent this (see \cite[Section 3]{ks} for proofs and details).

\begin{theorem}
\label{liechtenstein}
\textit{Let $G$ be a finite group, $N \unlhd G$ and $c$ be a $G$-stable block of $kN$. A $(c,G)$-Brauer pair is a pair $(Q, e_Q)$, where $Q \leq_p G$ with $\br^N_Q(c) \neq 0$.\\
We define the generalised Brauer category $\mathcal F_{(S,e'_S)}(G,N,c)$, where $(S,e'_S)$ is a maximal $(c,G)$-Brauer pair. Then $\mathcal F_{(S,e'_S)}(G,N,c)$ defines a fusion system on $S$. If $b$ is a block of $kG$ covering $c$ with maximal $b$-Brauer pair $(P,e_P)$, we have the relations $P \leq S$, $\mathcal F_{(P,e_P)}(G,b) \leq \mathcal F_{(S,e'_S)}(G,N,c)$ as well as $S \cap N=P \cap N$ and $\mathcal F_{(S \cap N,e'_{S \cap N})}(N,c) \dot{\unlhd} \mathcal F_{(S,e'_S)}(G,N,c)$}.
\end{theorem}

It should be noted, that we believe that the generalised Brauer category is very likely to appear in prospective stronger reduction theorems.

\begin{lemma}
\label{pizza}
\cite[Lemma 6.1]{ks} \textit{Let $G$ be a finite group with $N \unlhd G$ and $b$ be a block of $kG$ with defect group $D$. Then there exists a block $c$ of $kN$, which is covered by $b$, having $D \cap N$ as a defect group}.
\end{lemma}

We finish this section by making a link between fusion systems and blocks in stating the two known reduction theorems with respect to Conjecture \ref{seitan}.

\begin{theorem}
\label{radha}
\cite[Theorem 3.1]{solom} \textit{Let $\mathcal F$ be a reduction simple fusion system on a $p$-group $P$. Assume that $\aut(P)$ is a $p$-group. If $G$ is a finite group having an $\mathcal F$-block, then there exists a quasisimple group $L$ with $p'$-centre also having an $\mathcal F$-block}.
\end{theorem}

\begin{theorem}
\label{radhastancu}
\cite[Theorem 4.2]{ks} \textit{Let $\mathcal F_{1}$ and $\mathcal F_{2}$ be fusion systems on a $p$-group $P$ such that $\mathcal F_{1} \subseteq \mathcal F_{2}$. Assume that\\ 
(a) $\mathcal F_1$ is reduction simple,\\
(b) if $\mathcal F$ is a fusion system on $P$ containing $\mathcal F_{1}$, then $\mathcal F=\mathcal F_{1}$ or $\mathcal F=\mathcal F_{2}$,\\
(c) if $\mathcal F$ is a non-trivial normal subsystem of $\mathcal F_{2}$, then $\mathcal F=\mathcal F_{1}$ or $\mathcal F=\mathcal F_{2}$.\\
If there exists a finite group with an $\mathcal F_{1}$ or $\mathcal F_{2}$-block, then there also exists a quasisimple group $L$ with $p'$-centre with an $\mathcal F_{1}$ or $\mathcal F_{2}$-block}.
\end{theorem}

Note that both of these reductions can not be applied to the exotic Parker--Semeraro system $\mathcal{F}_7^1(6)$ (see \cite[Notation 5.14]{parkersem}): Since its subsystems are in correspondence to the subgroups of $C_6$, it will not be possible to fit each subsystem into a pair fulfilling the assumptions of Theorem \ref{radhastancu}. Since $\aut(S)$ is not a $7$-group, Theorem \ref{radha} is also not applicable.\\
However, in the next section, we generalise Theorem \ref{radhastancu} from $2$ fusion systems to all subsystems of $p$-prime index of a certain fusion system.

\section{Reduction Theorems}
We kick off this section by stating several reduction theorems which are essential for the study of fusion systems of blocks. These results are called respectively the First and Second Fong Reduction.

\begin{prop}
\label{fongone}
\cite[Part IV, Proposition 6.3]{ako} \textit{Let $\mathcal F$ be a fusion system on a $p$-group $P$ and let $G$ be a finite group having an $\mathcal F$-block $b$. Let $N \unlhd G$ and $c$ be a block of $kN$ which is covered by $b$. Then the group $I(c)=\{g \in G \mid ^gc=c \}$ has an $\mathcal F$-block}.
\end{prop}

We use the First Fong Reduction often in the following form:

\begin{cor}
\label{salted}
\textit{Let $\mathcal F$ be a fusion system and $G$ be a finite group possessing an $\mathcal F$-block $b$ such that $|G:Z(G)|$ is minimal among all finite groups having an $\mathcal F$-block. Then $b$ is inertial, i.e. it covers only $G$-stable blocks}.
\end{cor}

\textit{Proof}. Choose $N, c$ as in Proposition \ref{fongone}. Since $Z(G) \subseteq I(c)$, this proposition implies directly that $b$ is inertial. \hfill$\square$

\begin{theorem}
\label{fongtu}
\cite[Part IV, Theorem 6.6]{ako} \textit{Let $G$ be a finite group with $N \unlhd G$ and $c$ be a $G$-stable block of $N$ with trivial defect. Let $b$ be the block of $G$ covering $c$ and let $(P,e_P)$ be a maximal $b$-Brauer pair, then there exists a central extension $1 \rightarrow Z \rightarrow \widetilde{G} \rightarrow G/N \rightarrow 1$, where $Z$ is a cyclic $p'$-group. Furthermore, there is a block $\widetilde b$ of $k \widetilde G$ such that if we identify $P$ with the Sylow $p$-subgroup of the inverse image of $PN/N$ in $\widetilde G$, then there is a maximal $\widetilde b$-Brauer pair $(P,f_P)$ such that $\mathcal F_{(P,e_P)}(G,b)=\mathcal F_{(P,f_P)}(\widetilde G, \widetilde b)$}.
\end{theorem}

\begin{lemma}
\label{wiedersehen}
\textit{Let $\mathcal F$ be a reduction simple fusion system and $G$ be a finite group having an $\mathcal F$-block $b$ with non-abelian defect group $P$. If $G=\langle ^gP : g \in G \rangle$, then there exists a quasisimple group $L$ with $p'$-centre having an $\mathcal F$-block}.
\end{lemma}

\textit{Proof}.  We claim that if $N \unlhd G$ is proper, then $N$ has a block $d$ which is covered by $b$ and of defect zero. Indeed, by Lemma \ref{pizza}, we can choose $d$ such that it has $P \cap N$ as defect group. Since $N$ is normal and each morphism in $\mathcal F$ is induced by conjugation with an element in $G$, $P \cap N$ is also strongly $\mathcal F$-closed. If $P \cap N \neq 1$, then by reduction simplicity $P \cap N = P$, which means $N=G$ by assumption. This contradiction implies $P \cap N=1$.\\
By Theorem \ref{fongtu}, in this case there is a $p'$-central extension $\widetilde{G}$ of $G/N$ coming from an exact sequence $1 \rightarrow Z \rightarrow \widetilde{G} \rightarrow G/N \rightarrow 1$ having a block $c$ that is an $\mathcal{F}$-block. We now construct a quasisimple group $L$ with an $\mathcal{F}$-block. If we choose $N$ to be a maximal normal subgroup, then $G/N$ is either cyclic of prime order or $G/N$ is a non-abelian simple group. Note that by simplicity of $G/N$, we necessarily have $Z=Z(\widetilde G)$ in the extension above.\\
First, we assume that the first case holds, thus let $g \in G/N$ be the generating element and $\tilde{g}$ be a preimage of $g$ in $\widetilde{G}$. Then $\widetilde{G}=\langle Z, \tilde{g} \rangle$. This means $\widetilde{G}$ is abelian, hence so is $P$, which is a contradiction.\\
So, we are left with the case that $G/N$ is non-abelian simple. Define $L=[\widetilde{G},\widetilde{G}]$. We have $LZ/Z \unlhd \widetilde{G}/Z=G/N$. First assume $LZ/Z=1$, then $LZ=Z$, which means that $L \subseteq Z$, which implies $\widetilde{G}/Z \subseteq \widetilde{G}/L$. This is a contraction since $\widetilde{G}/L$ is abelian, but $G/Z$, which is a homomorphic image of $\widetilde{G}/L$, is not. So, by simplicity, $LZ/Z=\widetilde{G}/Z$, so $LZ=\widetilde{G}$. Taking commutators of this equation implies $[L,L]=[\widetilde{G},\widetilde{G}]=L$, so $L$ is perfect. Since we have $\widetilde{G}/Z=LZ/Z=\bigslant{L}{L \cap Z}$, $L$ is also a $p'$-central extension (of $G/N$) and thus quasisimple. \\
We have $\widetilde G=LZ=\bigslant{L \times Z}{K}$ for some $p'$-central subgroup $K$. This means there exists a fusion system preserving bijection between the blocks of $\widetilde G$ and the blocks of $L \times Z$ having $K$ in their kernel. But we can identify each block of $L \times Z$ with a block of $L$, since the blocks of $Z$ are just linear characters. It is easy to see that blocks which differ by a linear character give rise to the same fusion system. In particular, there is a block of $L$ which is an $\mathcal F$-block.\hfill$\square$\\

\textit{Proof of Theorem $\ref{youwereneverreallyhere}$}. Assume $G$ to be of minimal order among the groups possessing an $\mathcal F_H$-block $b$ for some $H \leq \Gamma_{p'}(\mathcal F)$. Let $N \unlhd G$ and $c$ be a block of $N$ covered by $b$. By Proposition \ref{fongone}, and our assumption, we may assume that $G=I(c)$. In particular, we may assume that $c$ is $G$-stable and the unique block of $N$ covered by $b$.\\
Now $P$ is a $(b,G)$-defect group. Consider $M:=\langle ^gP \mid g \in G \rangle \unlhd G$. Let $d$ be the block of $kM$ covered by $b$. Since $d$ is $G$-stable, we have a map $G \rightarrow \aut(kMd), g \mapsto c_g$, inducing a map $G \rightarrow \out(kMd)$. Let $K$ be the kernel of this map. Clearly, $M \subseteq K$. We claim $K=G$.\\
Indeed, let $f$ be the block of $kK$ covered by $b$. Let $(P,e_P)$ be a maximal $b$-Brauer pair and $(S,e'_S)$ a maximal $(G,K,f)$-Brauer pair such that $P \leq S$. By \cite[Section 5]{kuels}, $G/K$ is a $p'$-group. Thus, we may assume $P=S$. Furthermore $\mathcal F_H \leq \mathcal F_{(P,e'_P)}(G,K,f)$ and $\mathcal{F}_{(P,e'_P)}(K,f) \dot{\unlhd} \mathcal F_{(P,e'_P)}(G,K,f)$ by Theorem \ref{liechtenstein}.\\
Thus, by assumption $(b)$, $\mathcal F_{(P,e'_P)}(G,K,f)$ is of $p'$-index in $\mathcal F$ too and hence there is some $H' \leq \Gamma_{p'}(\mathcal F)$ such that $\mathcal F_{(P,e'_P)}(G,K,f)=\mathcal F_{H'}$. Similarly, by assumption $(c)$, there is also an $J \leq \Gamma_{p'}(\mathcal F)$ such that $\mathcal{F}_{(P,e'_P)}(K,f)=\mathcal F_J$. By the minimality of $G$, we deduce $G=K$.\\
By this observation, $G$ acts as inner automorphisms on $kMd$. Thus by \cite[Theorem 7]{hammer}, $kMd$ and $kGb$ have isomorphic source algebras. By \cite[Proposition 2.12]{solom}, we have that $d$ is an $\mathcal{F}_H$-block as well. Using the minimality once more, we obtain $G=M$. The previous lemma implies the theorem, since reduction simplicity of $\mathcal O^{p'}(\mathcal F)$ implies reduction simplicity of $\mathcal F$ for a fusion system $\mathcal F$. \hfill$\square$\\

Note that we obtain \cite[Theorem 4.2]{ks} as a corollary of the theorem we just proved by setting $\Gamma_{p'}(\mathcal F) = C_2$ respectively $\Gamma_{p'}(\mathcal F) = 1$.\\

We finish this section with proving Theorem \ref{cheeze}, which further restricts the structure of reduction simple fusion systems:\\

\textit{Proof of Theorem $\ref{cheeze}$}.  Let $\widetilde G$ be a group having an $\mathcal F$-block $\widetilde b$ subject to $|\widetilde G:Z(\widetilde G)|$ being minimal and $M \unlhd \widetilde G$ be maximal such that $P \not\subseteq M$. Note that since each normal $p$-subgroup is contained in each defect group of a $p$-block, we have that $O_p(\widetilde G) \leq P$. Furthermore, $O_p(\widetilde G)$ is strongly $\mathcal F$-closed, so either $O_p(\widetilde G)=1$ or $O_p(\widetilde G)=P$. If $O_p(\widetilde G)=P$, then $Z(P) \unlhd \widetilde G$. In particular, $Z(P)$ is strongly $\mathcal F$-closed, which is not possible since $P$ is non-abelian and $\mathcal F$ is reduction simple. Thus, $O_p(\widetilde G)=1$. In particular $Z(\widetilde G)$ is a $p'$-group. By maximality, we must have $Z(\widetilde G) \subseteq M$. Note that $P \cap M=1$ by reduction simplicity, but $P \cap M$ is a defect group of a block of $M$, which is covered by $\widetilde b$ by Lemma \ref{pizza}. In particular, there is a central extension $1 \rightarrow Z \rightarrow G \xrightarrow{\pi} \widetilde G/M \rightarrow 1$ for some central $p'$-group $Z$, such that $G$ has an $\mathcal F$-block $b$ by Theorem \ref{fongtu} (with the roles of $G$ and $\widetilde G$ interchanged).\\
For the first claim about $G$, note that $|G : Z(G)| \leq |G:Z|=|\widetilde G:M| \leq |\widetilde G:Z(\widetilde G)|$. In particular, by Corollary \ref{salted}, $b$ is inertial.\\
Suppose $H \unlhd G$ with $P \not\subseteq H$. We show $H \subseteq Z(G)$. We may assume $Z \subsetneq H$. Let $\varepsilon: \widetilde G \twoheadrightarrow \widetilde G/M$ be the canonical surjection. The maps $\varepsilon$ and $\pi$ induce bijections between the set of subgroups of $\widetilde G$ containing $M$ and the set of subgroups of $G$ containing $Z$, which preserves normality. This bijection sends $H$ to $\pi^{-1}(\varepsilon(H))$. Furthermore, $P \subseteq H$ if and only if $P \subseteq \pi^{-1}(\varepsilon(H))$. Since there is no normal subgroup of $\widetilde G$ properly containing $M$ and not containing $P$, it follows that there is no normal subgroup of $G$ properly containing $Z$ and not containing $P$, which implies our claim. In particular, $Z(G)$ is a $p'$-group.\\
Note that we have $O_p(G)=1$ for any $G$ having an $\mathcal F$-block. Since $F(G)=\prod \limits_{q \in \mathbb P} \prod \limits_{Q \in \syl_q(F(G))} Q$, we thus have $\syl_p(F(G))$ $=1$. Thus, by the above, $F(G) \subseteq Z(G)$, so in fact $F(G)=Z(G)$.\\
Now let $c$ be the block of $E(G)$ which is covered by $b$. If $E(G) \cap P=1$, $E(G) \subseteq Z(G)$ again by the above. But then $E(G)=1$ and $F(G)$ is central, which would mean that $G$ is abelian. So, $E(G) \cap P \neq 1$ and thus $P \subseteq E(G)=L_1 \cdots L_t$, where $\{L_1, \dots, L_t\}$ are the components of $G$. We have $E(G) \cong \bigslant{(L_1 \times \dots \times L_t)}{K}$ for $K \subseteq Z(L_1 \times \dots \times L_t)=Z(L_1) \times \dots \times Z(L_t)$. We claim that $K$ is a $p'$-group. It suffices to prove $O_p(L_i)=1$ for each $1 \leq i \leq t$. Indeed, if we assume the contrary, then the group $O_p(L_1) \cdots O_p(L_t)$ is a non-trivial normal subgroup of $E(G)$. In particular, $O_p(E(G)) \neq 1$. However, this is a characteristic subgroup of the layer, which implies $O_p(G) \neq 1$, a contradiction.\\
Thus, there is a fusion system preserving bijection between the blocks of $E(G)$ and the blocks of $L_1 \times \dots \times L_t$ which have $K$ in their kernel. In particular, we may assume $P = P_1 \times \dots \times P_t$ and $c = c_1 \times \dots \times c_t$, where for $1 \leq i \leq t$, $P_i$ is a defect group of the block $c_i$, which is a block of $L_i$ covered by $c$. If $r$ is the rank of $Z(P)$, then at most $r$ of the blocks $c_i$ can have non-trivial defect. Let $s \leq r$ be such that $s$ of the blocks $c_i$ have non-trivial defect. After possibly reordering, we may assume these are $c_{s+1}, \cdots, c_{t}$. We claim $L_{s+1} \cdots L_t \unlhd G$. Indeed, the conjugation action of $G$ on its components induces a group homomorphism $\sigma: G \rightarrow \sym(\{L_1, \dots, L_t\}) \cong \mathfrak S_t$ as follows: $\sigma(x)(i):=j$ iff $^xL_i=L_j$. Assume there is an $x \in G$ such that $^xL_i=L_j$ for $i \leq s$, $j>s$. Since $b$ is inertial, $c$ is $G$-stable. This means $^x c= c$, so $^x c_1 \times \dots \times$ $^x c_t = c_1 \times \dots \times c_t$, but this implies $^x P_i$ is non-trivial. This contradiction implies normality. Now we can apply the above to deduce $L_{s+1} \cdots L_t = 1$, which implies the third claim about $G$. \hfill$\square$\\

We can further restrict the structure of reduction simple fusion systems by specialising to the case of $Z(P)$ being cyclic:

\begin{theorem}
\label{glastonbury}
\textit{Let $P$ be a non-abelian $p$-group such that $Z(P)$ is cyclic and let $\mathcal F$ be a reduction simple fusion system on $P$. If $\mathcal F$ is block-realisable, then there exists a fusion system $\mathcal F_0$ on $P$ and a quasisimple group $L$ with an $\mathcal F_0$-block, where $\mathcal O_p(\mathcal F_0)=1$}.
\end{theorem}

\textit{Proof}. Assume $G$ is a finite group having an $\mathcal F$-block $b$ with defect group $P$. We may choose $G$ such that the conclusions of Theorem \ref{cheeze} hold. Let $L=\langle ^gP \mid g \in G \rangle$. Thus, since $P \subseteq E(G)$ as in the proof of Theorem \ref{cheeze}, we have $L \unlhd E(G)$. By Theorem \ref{cheeze}, the number of components of $G$ is bounded by the rank of $Z(P)$. By cyclicity of that group, $E(G)$ is quasisimple. Furthermore, $L$ is non-central, so we must have $L=E(G)$ is quasisimple.\\
Let $d$ be the block of $kL$ which is covered by $b$. Define $K$ to be the kernel of the map $G \rightarrow \out(kLd)$, which is induced by $G \rightarrow \aut(kLd), g \mapsto c_g$, and assume $K$ has a block $c$ which is covered by $b$. We get the triangle relations $\mathcal F \subseteq \mathcal F_{(P,e_P')}(G,K,c)$ and $\mathcal F_{(P,e_P')}(K,c) \dot{\unlhd} \mathcal F_{(P,e_P')}(G,K,c)=: \widetilde {\mathcal F}$ as in the proof of Theorem \ref{youwereneverreallyhere}. In the same fashion as in this theorem, application of \cite[Theorem 7]{kuels} and \cite[Proposition 2.12]{solom} also implies $\mathcal F_{(P,e_P')}(K,c) \cong \mathcal F_{(P,f_P)}(L,d)=:\mathcal F_0$. We have $\mathcal O_p(\mathcal F)=1$ by reduction simplicity. Thus, Lemma \ref{pint} implies $\mathcal O_p(\widetilde{\mathcal F})=1$ and Lemma \ref{ofscience} implies $\mathcal O_p(\mathcal F_0)=1$. \hfill$\square$

\section{The Parker--Semeraro Exotic Fusion Systems}
In this section, we use what we have developed so far to prove block-exoticity of the exotic Parker--Semeraro fusion systems.\\
In this chapter, $S$ will denote a Sylow $7$-subgroup of $G_2(7)$. As in the introduction, by a Parker--Semeraro system, we mean a fusion system $\mathcal F$ on $S$ such that $\mathcal O_p(\mathcal F) \neq 1$. In what follows, we will freely use the notation of \cite{parkersem}, where these systems have been classified.

We start by proving that it is not possible for most finite quasisimple groups to have a block with $S \in \syl_7(G_2(7))$ as a defect group. Fix $M$ to be the Monster group for the rest of this section.

\begin{prop}
\label{glaswegian}
\textit{Let $S \in \syl_7(G_2(7))$. Assume $G$ is a finite quasisimple group having a block with defect group $S$. Then either $G=M$ or $G=G_2(7)$}.
\end{prop}

Most work will have to be done to deal with groups of Lie type. We are going to restate a lemma from \cite{ks}, which will be very useful to deal with these groups.

\begin{lemma}
\label{handlotion}
\cite[Lemma 6.2]{ks} \textit{Let $H=LD$ be a finite group such that $L \unlhd H$ and $D$ is a $p$-group. Furthermore, let $c$ be a $D$-stable block of $kL$ with defect group $D \cap L$ and $\br^H_D(c) \neq 0$ and let $D'$ be such that $D \cap L \leq D' \leq D$. Then,\\
$(a)$ $c$ is a block of $LD'$ with defect group $D'$,\\
$(b)$ if $D'$ acts on $L$ as elements of $\inn(L)$, then $D'=(D' \cap L)C_{D'}(L)$}.
\end{lemma}

\begin{prop}
\textit{Let $G$ be a quasisimple finite group and denote the quotient $G/Z(G)$ by $\overline G$. Suppose $\overline{G}=G(q)$ is a finite group of Lie type and let $p$ be a prime number $\geq 7$, $p \neq q$. Let $D$ be a $p$-group such that $Z(D)$ is cyclic of order $p$ and $Z(D) \subseteq [D,D]$. If $D$ is a defect group of a block of $G$, then there are $n, k \in \mathbb{N}$ and a finite group $H$ with $\slg_n(q^k) \leq H \leq \gl_n(q^k)$ $($or $\su_n(q^k) \leq H \leq \gu_n(q^k))$ such that there is a block $c$ of $H$ with non-abelian defect group $D'$ such that $D'/Z$ is of order $|D/Z(D)|$ for some $Z \leq D' \cap Z(H)$}.
\end{prop}

\textit{Proof}. Suppose $G$ has a block with defect group $D$. Then the Sylow $p$-subgroups of $\overline{G}$ cannot be abelian, which implies that the Weyl group of the algebraic group corresponding to $\overline{G}$ has an order divisible by $p$, see \cite[Theorem 4.10.2(a)]{gls}. This implies that the exceptional part of the Schur multiplier of $\overline G$ is trivial, see \cite[Table 6.1.3]{gls}. Thus, there is a simple simply connected algebraic group $\overline K$ defined over $\overline{\mathbb{F}_q}$ and a Frobenius morphism $F: \overline{K} \rightarrow \overline{K}$ such that $\overline{K}^F$ is a central extension of $G$. If $\overline{K}$ is of type $A$, set $H:=\overline{K}^F$ and $c$ to be the block whose image has defect group $D$ under the algebra homomorphism $kH \rightarrow kG$ induced by $H \twoheadrightarrow G$.\\
Thus, assume $\overline{K}$ is not of type $A$. Since the kernel of $K^F \rightarrow \overline{K}^F$ is a $p'$-group, we have $|\overline K^F|_p=|K^F|_p$, so we may assume $G=\overline K^F$.\\
Denote the generator of $Z(D)$ by $z$. By Brauer's first main theorem, see \cite[Theorem 3.6]{radha}, we may assume the group $C_G(z)$ has a block $b$ with defect group $D$. Since $p$ is good, $\overline{K}$ is simply connected and thus $C_{\overline K}(G)$ is a Levi subgroup of $ \overline K$. If we denote $\overline Z:=Z(C_{\overline K}(z))^\circ$, it is a well-known fact that $C_{\overline K}(z)=[C_{\overline K}(z),C_{\overline K}(z)]\overline Z$. The latter commutator is simply connected and thus a direct product of its components, which are simply connected as well and permuted by $F$. In particular, we have a decomposition $[C_{\overline K}(z),C_{\overline K}(z)]=\prod \limits_{i=1}^{t} \prod \limits_{j=1}^{r_i} \overline{L}_{ij},$ where each $\overline{L}_{ij}$ is simply connected simple, and the set of these groups for a fixed $i$ lie in the same orbit.\\
Set $L_{i}=( \prod \limits_{j=1}^{r_i} \overline{L}_{ij})^F$. Then we have $C_G(z) = (L_1 \times \dots \times L_t)T^F$, for an $F$-stable maximal torus $T \leq C_{\overline K}(z)$.\\
Now $T^F$ is abelian and we have $\bigslant{D}{D \cap (L_1 \times \cdots \times L_t)} \cong \bigslant{T^F}{T^F \cap (L_1 \times \cdots \times L_t)}$. So $[D,D] \leq D \cap (L_1 \times \cdots \times L_t) \neq 1$. By Lemma \ref{pizza}, the latter group is furthermore a defect group of a block of $L_1 \times \dots \times L_t$. Now defect groups respect direct products and $Z(D) \cong C_p$. Thus we may assume $D \cap (L_1 \times \cdots \times L_t) \leq L_1$. In particular, $[D,D] \leq L_1$ and we may assume $Z(D) \leq L_1$. Since $Z(D)$ is central in $C_G(z)$, each $\overline{L}_{1j}$ is of type $A$ and Lie rank at least $p$, so $L_1$ is isomorphic to either $\slg_{n}(q^k)$ or $\su_{n}(q^k)$.\\
Let $x \in D \setminus Z(D)$. We want to show that $x$ does not centralise $L_1$. First note, if $\overline L=\slg_n(\overline{\mathbb{F}_q})$ and $\sigma$ a Frobenius morphism then $\overline{L}^\sigma$ is either $\slg_n(q^k)$ or $\su_n(q^k)$ for some $k$ and we have $C_{\overline L}(\overline L^\sigma) \leq Z(\overline L)$. Using the decomposition from above, we may write $x=(\prod \limits_{i=1}^{t} \prod \limits_{j=1}^{r_i} x_{ij})t_x$ for $x_{ij} \in \overline L_{ij}, t_x \in Z(\overline L)$. Let $y \in D$ such that $[x,y] \neq 1$ and write $y=(\prod \limits_{i=1}^{t} \prod \limits_{j=1}^{r_i} y_{ij})t_y$ for $y_{ij} \in \overline L_{ij}, t_y \in Z(\overline L)$. We have $[D,D] \leq L_1$, which means $[x_{11},y_{11}] \neq 1$, so $x_{11}$ does not centralise $\overline L_{11}$ as well as $\overline L_{11}^{F^{r_1}}$ by the above and thus also not $L_1$.\\
Now let $c$ be a block of $kL_1$ covered by $b$, then we may assume that $c$ is $D$-stable, $\br_D^{L_1D}(c) \neq 0$ and $D \cap L_1$ is a defect group of $c$, see Lemma \ref{pizza}. Let $D_0$ be the kernel of the map $D \rightarrow \out(L_1)$. Then $Z(D) \leq (D \cap L_1) \leq D_0$. Now if we apply Lemma \ref{handlotion}, we obtain $D_0=(D_0 \cap L_1)C_{D_0}(L_1)$. But we have seen that $C_D(L_1)=Z(D) \leq L_1$. So $D_0 \leq L_1$.\\
If $D_0=D$, we can take $H=L_1, D'=D$ and $Z=1$ and the claim holds. Thus, assume $D \neq D_0$. The elements of $T^F$ induce diagonal automorphisms of $L_1$. For special linear or unitary groups, the group of diagonal automorphisms modulo inner automorphisms is cyclic. In particular, $D/D_0$ is cyclic. Let $D/D_0=\langle \overline y \rangle$. Let $\eta \in \gl_n(q^k)_p$ (respectively $\gu_n(q^k)_p$) such that $^\eta u=$ $^y u$ for $u \in L_1$. In particular, $\eta$ stabilises $c$. Define $H:=L_1 \langle \eta \rangle$ to obtain $H$ as in the claim. Furthermore, define $D':=\langle D_0, \eta \rangle \leq H$. We also have $C_{L_1}(D)=C_{L_1}(D')$, so $\br_{D'}^H(c) \neq 0$. Now $H=L_1 D'$ and $D_0 (\cong D' \cap L_1$ by construction) is a defect group of $c$ as a block of $kL_1$. Thus, we can apply Lemma \ref{handlotion} to obtain that $c$ is a block of $kH$ with defect group $D'$.\\
We have $D'=\langle D_0, \eta \rangle $ and $D=\langle D_0, y \rangle$. The canonical maps $D' \rightarrow \aut(L_1)$ and $D \rightarrow \aut(L_1)$ have the same image. Thus, $\bigslant{D'}{C_{D'}(L_1)} \cong \bigslant{D}{C_{D}(L_1)} = D/Z(D)$. Define $Z:=Z(\gl_n(q^k)) \cap D' = C_{D'}(L_1)$, which gives $D'/Z  \cong D/Z(D)$. In particular, $D'$ is non-abelian since $y$ acts non-trivially on $D_0$.
\hfill$\square$

\begin{prop}
\label{hjelp}
\textit{If $G$ is as in the previous proposition, then $G$ has no blocks with defect groups isomorphic to a Sylow $7$-subgroup of $G_2(7)$}.
\end{prop}

\textit{Proof}. Assume $D \in \syl_7(G_2(7))$, in particular $|D|=7^6$ and $Z(D) \cong C_7$. Let $H$, $D'$ be as in the assertion of the previous proposition with $p=7$. Assume first $H \leq \gl_n(q^k)$ and let $a$ be such that $|q^k-1|_7=7^a$. Then, since $\slg_n(q^k) \leq H$, we have $|D'|=|D/Z(D)| \cdot |Z|=7^5|Z| \leq 7^5|Z(H)| \leq 7^5|Z(\slg_n(q^k))|\leq 7^{5+a}$. Now the block of $k\gl_n(q^k)$ covering $c$ has a defect group of order at most $7^{2a+5}$. But it is a well-known fact, that (non-abelian) defect groups of $\gl_n(q^k)$ have order at least $7^{7a+1}$. Thus, $7^{7a+1} \leq 7^{2a+5}$, which is a contradiction. The case $H \leq \gu_n(q^k)$ can be shown in the same fashion by considering the $7$-part of $q^k+1$ instead of $q^k-1$. \hfill$\square$\\

We use this observation to prove Proposition \ref{glaswegian}.\\

\textit{Proof of Proposition $\ref{glaswegian}$}. By the previous proposition, we may assume $G$ is not a group of Lie type over characteristic not equal to $7$.\\
First, assume $G/Z(G)$ is an alternating group $\mathfrak A_m$. Then $S$ is isomorphic to a Sylow $7$-subgroup of some symmetric group $\mathfrak S_{7w}$ with $w \leq 6$. Define the cycle $\sigma_i=((i-1)7+1,\dots,i7)$ and the subgroup $S'=\langle \sigma_1, \dots, \sigma_6 \rangle\leq \mathfrak A_m$. Then $S' \in \syl_7(\mathfrak A_m)$. But this group is abelian, which means that $S \notin \syl_7(\mathfrak A_m)$.\\
Next, assume $G$ is a group of Lie type over a field of characteristic $p=7$. In this case, $Z(G)$ is a $7'$-group and we may assume that $G=\mathbf G^F(q)$, where $\mathbf G$ is simple and simply connected, $F$ is a Frobenius morphism and $q=p^f$ for some $f \in \mathbb{N}$. Furthermore, $S \in \syl_7(G)$ by \cite[Theorem 6.18]{caen}. We first deal with the classical groups. First consider type $B_n$, here we have $|G|_p=q^{n^2}$, which can be equal to $7^6$ only if $n=1$, which is not possible since $n>1$ for these groups. The case is the same for the groups of type $C_n$. For type $D_n$ and $^2 D_n$, we have $|G|_p=q^{n(n-1)}$. Since $n>3$, $p^6$ is also no possibility here. Finally, consider types $A_n$ and $^2 A_n$, here we have $|G|_p=q^{\frac{1}{2}(n+1)n}$. If $n=1$, these groups have abelian Sylow subgroups and if $n \geq 4$, the order is too big. Thus, the only possibilities are $n=2$ or $n=3$. However, in these cases we obtain Sylow $7$-subgroups which are conjugate to the groups of upper unitary triangle matrices. Hence, these groups have nilpotency class $2$ respectively $3$. However, the nilpotency class of $S$ is $5$. This leaves us with the exceptional groups of Lie type. Looking at their orders, we can directly exclude the exceptional Steinberg groups, the Suzuki, Ree and Tits groups. For the exceptional Chevalley groups, $G_2(7)$ is the only possibility.\\
Finally, if $G/Z(G)$ is sporadic, the monster $M$ is the only group with a $7$-part big enough to contain $S$, which implies our claim. \hfill$\square$\\

This result can be used to achieve a reduction specifically for the Parker--Semeraro systems.

\begin{lemma}
\label{nauru}
\textit{Let $\mathcal F$ be an exotic Parker--Semeraro system. Then $\mathcal F$ is reduction simple}.
\end{lemma}

\textit{Proof}. Assume $1 \neq N \leq S$ is strongly $\mathcal F$-closed. In particular $N \unlhd S$, which implies $Z(S) \leq N$. Thus, as in the proof of \cite[Theorem 6.2]{parkersem}, we obtain $N=S$. \hfill$\square$

\begin{theorem}
\textit{Suppose there is an exotic Parker--Semeraro system $\mathcal F$ which is block-realisable. Then there is an exotic Parker--Semeraro system $\mathcal F_0$ which is block-realisable by a block of a quasisimple group}.
\end{theorem}

\textit{Proof}. Assume $G$ is a finite group having an $\mathcal F$-block $b$. We may choose $G$ such that the conclusions of Theorem \ref{cheeze} hold. Define $L = \langle ^g S \mid g \in G \rangle \unlhd G$. Since $Z(S)$ is cyclic of order $p$, $\mathcal F$ satisfies the hypothesis of Theorem \ref{glastonbury} with $P=S$. Arguing as in the proof of that theorem, $L=E(G)$ is quasisimple, and there is a block of $kL$ with defect group $S$ and fusion system $\mathcal F_0$ such that $O_p(\mathcal F_0)=1$. So, $\mathcal F_0$ is a Parker--Semeraro system.\\
If $\mathcal F_0$ is exotic, we are done. Suppose $\mathcal F_0$ is realisable. So, either $\mathcal F_0=\mathcal F_S(M)$ or $\mathcal F_0= \mathcal F_S(G_2(7))$. By Proposition \ref{glaswegian}, $L=M$ or $L=G_2(7)$ and hence $L$ is simple.\\
We claim $G=L$. Indeed, consider the map $\varphi: G \rightarrow \out(F^\ast(G))$ and let $g \in \ker(\varphi)$. Then there exists $x \in F^\ast(G)$ such that $y \in F^\ast(G)$ with $gyg^{-1}=xyx^{-1}$, i.e. $x^{-1}g \in C_G(F^\ast(G))=Z(F^\ast(G))$. This implies $\ker \varphi=F^\ast(G)$ and thus $G/F^\ast(G) \cong \out_G(F^\ast(G))$. But in our case, we have $F^\ast(G)=Z(G)E(G)$, so $\out_G(F^\ast(G))=\out_G(E(G))$. However, this group needs to be trivial since $\out(M)=\out(G_2(7))=1$. This implies $G=Z(G)E(G)$ and in either case $G=Z(G) \times L$. By choice of $G$, the claim follows.\\
We now claim that in either of these cases, $b$ needs to be the principal block. For $G_2(7)$ this is a well known fact, see e.g. \cite[Example 3.8]{radhagunter}. So, assume $G=M$. We want to compute $C_{N_M(P)}(S)$, where $P$ is a subgroup of order $7$ of the monster. We know that $C_{N_M(P)}(S)=Z(S) \times O_{7'}(C_{N_M(P)}(S))$. Furthermore, by \cite[Theorem 1.1]{weeddorn}, $N_M(P)$ is $7$-constrained, i.e. $C_{N_M(P)}(O_7(N_M(P))) \subseteq O_7(N_M(P))$. Now $N_M(P)=7^{1+4}(2 \mathfrak{S}_7 \times 3)$, i.e. $O_7(N_M(P))=7^{1+4}$ and thus $C_{N_M(P)}(S) \cong C_p$. But this means $C_M(S)/Z(S)$ is trivial and since the $N_M(S)$-classes of characters of this group are in 1:1-correspondence with blocks of $M$ with defect group $S$, the claim follows. However, this implies that $\mathcal F$ cannot be exotic, which is a contradiction. \hfill$\square$\\

We use this to deduce block-exoticity of the Parker--Semeraro systems.\\

\textit{Proof of Theorem $\ref{oh hai mark}$}. Let $S$ be a Sylow $7$-subgroup of $G_2(7)$ and let $\mathcal{F}$ be one of the exotic Parker--Semeraro systems. Assume $G$ is a group with an $\mathcal{F}$-block. By the previous theorem, we may assume that $G$ is quasisimple. Let $A$ be the simple quotient of $G$. By Proposition \ref{glaswegian}, we may assume either $A=G_2(7)$ or $A=M$, and thus since $\out(A)=1$ in both these cases, also either $G=G_2(7)$ or $G=M$. But as in the proof of the previous theorem, this means that the $\mathcal F$-block is principal, which is not possible for an exotic fusion system.
 \hfill$\square$

\end{document}